# Smoothing analysis of two-color distributive relaxation for solving 2D Stokes flow by multigrid method

Xingwen Zhu[a*], Lixiang Zhang[b]

[a] *School of Mathematics and Computer, Dali University, Dali 671003, Yunnan, China.*
[b] *Department of Mechanical Engineering, Kunming University of Science and Technology, Kunming 650500, Yunnan, China.*

**Abstract**

Smoothing properties of two-color distributive relaxation for solving a two-dimensional (2D) Stokes flow by multigrid method are theoretically investigated by using the local Fourier analysis (LFA) method. The governing equation of the 2D Stokes flow in consideration is discretized with the non-staggered grid and an added pressure stabilization term with stabilized parameters to be determined is introduced into the discretization system in order to enhance the smoothing effectiveness in the analysis. So, an important problem caused by the added pressure stabilization term is how to determine a suitable zone of parameters in the added term. To that end, theoretically, a two-color distributive relaxation, developed on the two-color Jacobi point relaxation, is established for the 2D Stokes flow. Firstly, a mathematical constitution based on the Fourier modes with various frequency components is constructed as a base of the two-color smoothing analysis, in which the related Fourier representation is presented by the form of two-color Jacobi point relaxation. Then, an optimal one-stage relaxation parameter and related smoothing factor for the two-color distributive relaxation are applied to the discretization system, and an analytical expression of the parameter zone on the added pressure stabilization term is established by LFA. The obtained analytical results show that numerical schemes for solving 2D Stokes flow by multigrid method on the two-color distributive relaxation have a specific convergence zone on the parameters of the added pressure stabilization term, and the property of convergence is independent of mesh size, but depends on the parameters of the pressure stabilization term.

## 1. Introduction

Multigrid methods [1-7] are generally considered as one of the fastest numerical methods which have an optimally computational complexity for solving partial differential equations (PDEs), especially, for solving 2D Stokes flow governed by the following equations,

$$\begin{cases} -\Delta \vec{u} + \vec{\nabla} p = \vec{f} & (x,y) \in \Omega \\ \vec{\nabla} \cdot \vec{u} = 0 & (x,y) \in \Omega \\ \vec{u} = \vec{g} & (x,y) \in \partial\Omega \end{cases} \quad (1)$$

where $(x,y) \in \Omega \subseteq \mathbb{R}^2$, $\vec{u} = (u(x,y), v(x,y))^T$ is the velocity field, $p = p(x,y)$ represents the pressure, $\vec{f} = (f_1(x,y), f_2(x,y))^T$ is the external force field, and $\partial\Omega$ is the boundary of the computational domain on which the Dirichlet boundary condition is assigned.

In multigrid methods, smoothing relaxations play an important role for numerical simulations. Several relaxation techniques have been developed for solving systems of PDEs, which are roughly classified into two categories, collective and decoupled relaxations [8]. The collective relaxation is considered as straightforward generalization of a scalar system [2], and the decoupled relaxation is based on a distributive Gauss-Seidel technique [9] and generalized to analyze incomplete LU factorization for Stokes and Navier-Stokes systems [10]. Recently, a distributive Gauss-Seidel relaxation based on the least squares commutator has been used for solving the Stokes flow system [11]. More relaxation techniques for solving Stokes flow are presented in [12] and [13].

* Corresponding author. Tel.: ++86 872 2219600; fax: ++86 872 2219681.
E-mail address: zxw4688@126.com.

In smoothing analysis, the local Fourier analysis (LFA) is a quite useful tool to design efficient algorithms and to predict convergence factors for solving PDEs with high order accuracy [1-7]. In [14], a distributive relaxation method improved by using LFA is used for solving governing equations of poroelasticity. A multigrid solver on LFA for the Navier-Stokes equations is designed in [15]. All-at-once multigrid approaches for optimal control problems with LFA are discussed in details in [16,17], and an analytical expression of the convergence factors is given by using symbolic computation in [18]. With the assistance of LFA, a general definition on the multicolor relaxation is provided in [19]. Smoothing analysis of the two-color relaxation with LFA is discussed in [20-25], and the four-color relaxation with tetrahedral grids is presented in [26, 27]. Analyses of the solutions to the Stokes flow are discussed in literature due to its wide range in engineering applications [1,2,6,8]. Stabilization for nonstaggered discretizations can be achieved by adding an artificial elliptic pressure term, $-ch^2\Delta_h p$, to the discretization of continuity equation $\vec{\nabla}\cdot\vec{u}=0$, which can also maintain second-order accuracy of discretizations in PDE system, here $\Delta_h$ is the standard central difference of Laplace operator $\Delta$. Excellent performances of the multigrid with LFA for solving Stokes flow were well known in past years, for example, in [1-8, 10, 13], where the specific parameter values of $c = 1/16$ and $c = 1/8$ for the convergence properties of the numerical results are discussed, and, the mathematical expression of $h$-ellipticity of the discrete Stokes equations on the parameter $c$ was given in [1]. However, so far the stabilization zone of the parameter $c$ and its mathematical expression have been still missing.

In this paper, we will investigate the relation of parameter $c$ and the smoothing factor of multigrid method for 2D Stokes equations by means of the two-color Jacobi distributive relaxation, and theoretically analyze the smoothing factor of the stabilization scheme for the 2D Stokes flow in the stabilized central difference discretization on a non-staggered grid. To that end, we first employ two different Fourier modes to construct a Fourier invariant subspace ($2h$-harmonics subspace) for using the two-color distributive Jacobi relaxation in LFA, through which a mathematical relation of the smoothing factor and the parameter zone of the relaxation parameter with $c > 0$ are obtained by virtue of the optimal one-stage relaxation method. Then based on the two-color distributive relaxation, we apply the smoothing factor and relaxation parameter to the studied 2D Stokes flow, and analyze the theoretical form of the parametric expression of the added pressure term by LFA. The obtained analytical results show that the two-color distributive relaxation technique works well for the smoothing analysis of the 2D Stokes flow in the framework of multigrid, where a specific convergence zone on the pressure stabilization parameters c is achieved for the first time. Moreover, we can also see that the property of convergence is independent of mesh size, but depends on the pressure stabilization parameter $c$.

## 2. Discretizing Stokes flow and LFA

*2.1 Discrete Stokes flow*

From (1), a 2D Stokes operator can be written as

$$L = \begin{pmatrix} -\Delta & 0 & \partial_x \\ 0 & -\Delta & \partial_y \\ \partial_x & \partial_y & 0 \end{pmatrix}. \tag{2}$$

On a non-staggered grid

$$G_h = \{\vec{x} = (x,y) := (k_1 h, k_2 h) \mid (k_1, k_2) \in \mathbb{Z}^2\}, \tag{3}$$

the discretization of Stokes operator (2) by means of the standard central difference stencil is given as

$$L'_h = \begin{pmatrix} -\Delta_h & 0 & \partial_x^h \\ 0 & -\Delta_h & \partial_y^h \\ \partial_x^h & \partial_y^h & 0 \end{pmatrix}, \tag{4}$$

where $h$ denotes the uniform mesh size of spatial discretization, $-\Delta_h$, $\partial_x^h$, and $\partial_y^h$ are the second-order difference operator with the following discrete stencils

$$-\Delta_h = \frac{1}{h^2}\begin{bmatrix} & -1 & \\ -1 & 4 & -1 \\ & -1 & \end{bmatrix}_h, \quad \partial_x^h = \frac{1}{2h}\begin{bmatrix} -1 & 0 & 1 \end{bmatrix}_h, \quad \partial_y^h = \frac{1}{2h}\begin{bmatrix} 1 \\ 0 \\ -1 \end{bmatrix}_h. \tag{5}$$

In fact, it is well known that the above non-staggered-grid scheme (4) is unstable [1,2,6], we need to add artificially elliptic pressure term, $-ch^2\Delta_h$, into the continuity equation (1) to improve the numerical stabilization in computation, where $c > 0$ is stabilization parameter. Thus, the above discretization operator (4) is changed to

$$L_h = \begin{pmatrix} -\Delta_h & 0 & \partial_x^h \\ 0 & -\Delta_h & \partial_y^h \\ \partial_x^h & \partial_y^h & -ch^2\Delta_h \end{pmatrix}. \tag{6}$$

*2.2 Elements of LFA in multigrid*

A crucial point of using the multigrid method is to identify multigrid components which are used to construct an efficient interplay between the relaxation and coarse grid corrections. A useful tool for a proper selection of the components is the local Fourier analysis (LFA). In [1-7, 14], LFA is applied to develop efficient multigrid methods for solving linear elliptic equations with constant (or frozen) coefficients , which is based on a simplification that boundary conditions are neglected, and all occurring operators are extended to an infinite grid. On an infinite grid in the numerical simulation, the approximation and corresponding error and residual are represented by linear combinations of certain exponential functions or by the Fourier mode functions, which are used as a unitary basis in the space of the bounded infinite grid functions [1-7].

Actually, on the non-staggered grid $G_h$ defined in (3), a unitary basis based on Fourier modes is defined as

4[1,2],

$$\varphi_h(\vec{\theta}, \vec{x}) := exp(i\vec{\theta} \cdot \vec{x} / h), \tag{7}$$

where $\vec{\theta} = (\theta_1, \theta_2) \in \Theta := (-\pi, \pi]^2$ is called the Fourier frequency, $\vec{x} \in G_h$, and $i = \sqrt{-1}$. The Fourier space is spanned as

$$F(\vec{\theta}) := span\{\varphi_h(\vec{\theta}, \vec{x}) \mid \vec{\theta} \in \Theta\}. \tag{8}$$

Based on the grid $G_h$ in (3) and apply (7), the Fourier symbol of a 2D scalar discrete operator $D_h$ that is equipped the following discrete stencil

$$D_h = [l_{\vec{k}}]_h, \tag{9}$$

where $l_{\vec{k}} \in \mathbb{R}$, and $\vec{k} \in J \subset \mathbb{Z}^2$ containing $(0,0)$, can be defined as [1-7]

$$\tilde{D}_h(\vec{\theta}) := \sum_{\vec{k} \in J} l_{\vec{k}} exp(i\vec{\theta} \cdot \vec{k}), \tag{10}$$

with $\vec{\theta} \cdot \vec{k} = \theta_1 k_1 + \theta_2 k_2$, subjecting to

$$D_h \varphi_h(\vec{\theta}, \vec{x}) = \tilde{D}_h(\vec{\theta}) \varphi_h(\vec{\theta}, \vec{x}). \tag{11}$$

In the LFA for a two-dimensional case of the standard coarsening [1,2], the low frequencies are given by $\vec{\theta} = (\theta_1, \theta_2) \in \Theta_{low} = (-\frac{\pi}{2}, \frac{\pi}{2}]^2$, then the high frequencies by $\vec{\theta} \in \Theta_{high} = \Theta \setminus \Theta_{low}$. An ideal coarse-grid correction operator $Q_h^{2h}$ is introduced in [1,2], where the low-frequency Fourier modes are dropped out while the high-frequency Fourier modes remain. Thus, a general coarsening strategy is stated follows

$$Q_h^{2h} \varphi_h(\vec{\theta}, \vec{x}) := \begin{cases} \varphi_h(\vec{\theta}, \vec{x}) & \vec{\theta} \in \Theta_{high} \\ 0 & \vec{\theta} \in \Theta_{low} \end{cases}, \tag{12}$$

In this paper, the standard coarsening in multigrid is applied implicitly. The main work of LFA is to analyze different multigrid relaxations by evaluating their effects on the Fourier modes.

**3. Fourier representation of two-color distributive relaxation**

To develop the Fourier representation of the two-color distributive relaxation by LFA, we need to divide the grid $G_h$ (3) into two disjoint subsets $G_h^R$ and $G_h^B$, referring to as the red and black points, respectively. Let $S_h$ denote a certain point relaxation, such as Jacobi point relaxation, Gauss-Seidel point relaxation, and etc. then, a two-step process [2] is required to construct a complete two-color relaxation $S_h^{RB}(\omega)$. In the first step, $S_h^R(\omega)$, the unknowns located at the red points are the only ones to be smoothed, whereas the unknowns at the black points remain to be unchanged. In the second step $S_h^B(\omega)$, the unknowns at the black points are changed by using the new values calculated at the red points in the first step. So, a complete red-black point process is obtained by the following iteration

$$S_h^{RB}(\omega) = S_h^B(\omega) S_h^R(\omega), \tag{13}$$

From the above process, it is noted that the Fourier modes (7) are no longer eigenfunctions of (13) on the grid $G_h$ (3). However, the relaxation operator (13) leaves certain low-dimensional subspaces of Fourier modes invariant, yields a block-diagonal Fourier representation of the smoothing operator that consists of small blocks, such that the smoothing factor can be calculated, easily. In the following descriptions, a smoothing process of the two-color point relaxation is presented by LFA.

In the process of the two–color point relaxation, the grid $G_h$ (3) is divided into two disjoint subsets $G_h^0$ and $G_h^1$, i.e. $G_h = G_h^0 \cup G_h^1$ with

$$G_h^\beta = \{\vec{x} = (k_1 h_1, k_2 h_2) \mid k_1 + k_2 = \beta \bmod 2, \vec{k} \in \mathbb{Z}^2\} \tag{14}$$

where $\beta = 0, 1$. According to (12), the Fourier space (8) is subdivided into the corresponding 2h-harmonics subspace

$$F_{2h}^2(\vec{\theta}) := span\left\{\varphi_h(\vec{\theta}^0, \vec{x}), \varphi_h(\vec{\theta}^1, \vec{x})\right\} \tag{15}$$

with

$$\vec{\theta}^\alpha = (\vec{\theta} + (\alpha, \alpha)\pi) \bmod 2\pi, \alpha = 0, 1, \vec{\theta} \in \Theta_{low}. \tag{16}$$

From (15), if the relaxation operator $S_h$ satisfies

$$S_h(\varphi_h(\vec{\theta}^0, \vec{x}), \varphi_h(\vec{\theta}^1, \vec{x})) = (\varphi_h(\vec{\theta}^0, \vec{x}), \varphi_h(\vec{\theta}^1, \vec{x})) \widehat{S}_h(\vec{\theta}) \tag{17}$$

i.e. $S_h : F_{2h}^2(\vec{\theta}) \to F_{2h}^2(\vec{\theta})$, the matrix $\widehat{S}_h(\vec{\theta})$ is called Fourier representation of $S_h$ on the grid $G_h$ (3). Next, the Fourier representation of (13) on the 2h-harmonics subspace (15) is discussed. First, one constitution among the various Fourier modes defined by (14)-(16) is presented as follows.

**Proposition 1** For $\forall \alpha, \beta \in \{0,1\}, \forall \vec{x} \in G_h^\beta$, if $\vec{\theta} \in \Theta_{low}$, then the following formulation holds

$$\varphi_h(\vec{\theta}^\alpha, \vec{x}) = exp(i\pi\alpha\beta)\varphi_h(\vec{\theta}, \vec{x}). \tag{18}$$

Proof. From (7) and (14), for $\forall \vec{x} \in G_h^\beta \subseteq G_h$, the following equation holds

$$\varphi_h(\vec{\theta}^\alpha, \vec{x}) = exp(i\vec{\theta}^\alpha \cdot \vec{k}), \tag{19}$$

where $\vec{k} = (k_1, k_2) \in \mathbb{Z}^2$. From (16), $\exists \vec{n} = (n_1, n_2) \in \mathbb{Z}^2$ subjects to

$$\vec{\theta}^\alpha \cdot \vec{k} = \vec{\theta} \cdot \vec{k} + \pi\alpha(k_1 + k_2) + 2\pi\vec{n} \cdot \vec{k}. \tag{20}$$

For $\vec{x} \in G_h^\beta$, from (14), $\exists p \in \mathbb{Z}$ subjects to $k_1 + k_2 = \beta + 2p$. Substitute (20) into (19) and apply (7), yield

$$\begin{aligned}\varphi_h(\vec{\theta}^\alpha,\vec{x}) &= exp[i(\vec{\theta}\cdot\vec{k}+\pi\alpha\beta+2\pi\alpha p+2\pi\vec{n}\cdot\vec{k})]\\ &= exp[i(\vec{\theta}\cdot\vec{k}+\pi\alpha\beta]\\ &= \varphi_h(\vec{\theta},\vec{x})exp(i\pi\alpha\beta)\end{aligned} \qquad (21)$$

The proposition is thus proved.□

Subsequently, the smoothing analysis process of the two-color relaxation on the subspace of the 2$h$-harmonics (15) can be carried out as follows. By (13) and (14) and without loss of generality, let $G_h^0$ and $G_h^1$ correspond to $G_h^R$ and $G_h^B$, respectively, thus (13) can be rewritten as

$$S_h^{01}(\omega) = S_h^1(\omega)S_h^0(\omega). \qquad (22)$$

Apply Proposition 1 and (17), the Fourier representation of (22) can be obtained as shown in the following theorem.

**Theorem 2** The iteration operator $S_h^{01}(\omega)$ for the two-color relaxation leaves the subspace of the 2$h$-harmonics (15) invariant.

Proof. From the process of the two-color relaxation, the operator $S_h^\beta(\omega)$ of on the grid $G_h$ (3) leads to

$$S_h^\beta(\omega)\varphi_h(\vec{\theta},\vec{x}) = \begin{cases} \tilde{S}_h^\beta(\vec{\theta},\omega)\varphi_h(\vec{\theta},\vec{x}) & \forall \vec{x}\in G_h^\beta,\\ \varphi_h(\vec{\theta},\vec{x}) & \forall \vec{x}\notin G_h^\beta, \end{cases} \qquad (23)$$

where $\tilde{S}_h^\beta(\vec{\theta},\omega)$ is the Fourier symbol of the point relaxation $S_h^\beta(\omega)$ on the grid $G_h^\beta$ (14) with $\beta=0,1$.

From (15) and (17), we need to find out two complex numbers, $a_0$ and $a_1$, $\forall\,\alpha,\beta\in\{0,1\}$, subjecting to

$$S_h^\beta(\omega)\varphi_h(\vec{\theta}^\alpha,\vec{x}) = a_0\varphi_h(\vec{\theta}^0,\vec{x})+a_1\varphi_h(\vec{\theta}^1,\vec{x}). \qquad (24)$$

Taking $A_\alpha^\beta = \tilde{S}_h^\beta(\vec{\theta}^\alpha,\omega)$, and applying (23), (24) and Proposition 1, we can present the following linear equations with respect to $a_0$ and $a_1$,

$$\begin{cases} a_0+a_1 exp(i\beta\pi) = A_\alpha^\beta exp(i\alpha\beta\pi), & \forall\vec{x}\in G_h^\beta,\\ a_0+a_1 exp[i(1-\beta)\pi] = exp[i\alpha(1-\beta)\pi], & \forall\vec{x}\notin G_h^\beta, \end{cases} \qquad (25)$$

where $\alpha,\beta\in\{0,1\}$. Apply (15), (17) and (24), and solve the system of linear equation for $a_0$ and $a_1$, the Fourier representations of iteration operators $S_h^0(\omega)$ and $S_h^1(\omega)$ can be obtained as

$$\widehat{S}_h^0(\vec{\theta},\omega) = \frac{1}{2}\begin{pmatrix} A_0^0+1 & A_1^0-1\\ A_0^0-1 & A_1^0+1 \end{pmatrix},\quad \widehat{S}_h^1(\vec{\theta},\omega) = \frac{1}{2}\begin{pmatrix} A_0^1+1 & -A_1^1+1\\ -A_0^1+1 & A_1^1+1 \end{pmatrix}, \qquad (26)$$

where $A_\alpha^\beta = \tilde{S}_h^\beta(\vec{\theta}^\alpha,\omega)$ and $\alpha,\beta\in\{0,1\}$. Furthermore, the Fourier representations of the two-color relaxation (22) is attained from (26) as follows,

$$\widehat{S}_h^{01}(\vec{\theta},\omega) = \widehat{S}_h^{1}(\vec{\theta},\omega)\widehat{S}_h^{0}(\vec{\theta},\omega) = \frac{1}{2}\begin{pmatrix} A_0^1+1 & -A_1^1+1 \\ -A_0^1+1 & A_1^1+1 \end{pmatrix} \cdot \frac{1}{2}\begin{pmatrix} A_0^0+1 & A_1^0-1 \\ A_0^0-1 & A_1^0+1 \end{pmatrix}. \quad (27)$$

Considering (17), we thus prove the theorem.□

## 4. Smoothing analysis of two-color distributive relaxation for Stokes flow

### 4.1. Distributive relaxation of discrete Stokes flow

A distributive relaxation for the discrete system (6) can be constructed as follows [1,2,7]. First, the distributive operator is given as

$$C_h = \begin{pmatrix} I_h & 0 & -\partial_x^h \\ 0 & I_h & -\partial_y^h \\ 0 & 0 & -\Delta_h \end{pmatrix}, \quad (28)$$

where $I_h$ is unit operator with discrete stencil $[1]_h$. Then, the discrete system (6) is transformed follows (28)

$$L_h C_h = \begin{pmatrix} -\Delta_h & 0 & 0 \\ 0 & -\Delta_h & 0 \\ \partial_x^h & \partial_y^h & ch^2 \Delta_h^2 - \Delta_{2h} \end{pmatrix}, \quad (29)$$

with the discrete stencils

$$\Delta_h^2 = \frac{1}{h^4}\begin{bmatrix} & & 1 & & \\ & 2 & -8 & 2 & \\ 1 & -8 & 20 & -8 & 1 \\ & 2 & -8 & 2 & \\ & & 1 & & \end{bmatrix}_h, \quad -\Delta_{2h} = \frac{1}{4h^2}\begin{bmatrix} 0 & 0 & -1 & 0 & 0 \\ 0 & 0 & 0 & 0 & 0 \\ -1 & 0 & 4 & 0 & -1 \\ 0 & 0 & 0 & 0 & 0 \\ 0 & 0 & -1 & 0 & 0 \end{bmatrix}_h = \frac{1}{4h^2}\begin{bmatrix} & -1 & \\ -1 & 4 & -1 \\ & -1 & \end{bmatrix}_{2h}. \quad (30)$$

Next (9)-(11) result in the following Fourier symbols of the scalar discrete operators of (30),

$$\tilde{\Delta}_h^2(\vec{\theta}) = (-\tilde{\Delta}_h(\vec{\theta}))^2, \quad -\tilde{\Delta}_{2h}(\vec{\theta}) = -[\tilde{\partial}_x^h(\vec{\theta})]^2 - [\tilde{\partial}_y^h(\vec{\theta})]^2, \quad (31)$$

where

$$-\tilde{\Delta}_h(\vec{\theta}) = \frac{1}{h^2}(4 - exp(-i\theta_1) - exp(i\theta_1) - exp(-i\theta_2) - exp(i\theta_2))$$
$$= \frac{1}{h^2}(4 - 2\cos\theta_1 - 2\cos\theta_2) \quad (32)$$

$$\tilde{\partial}_x^h(\vec{\theta}) = \frac{1}{2h}(exp(i\theta_1) - exp(-i\theta_1)) = \frac{1}{h}i\sin\theta_1 \quad (33)$$

$$\tilde{\partial}_y^h(\vec{\theta}) = \frac{1}{2h}(exp(i\theta_2) - exp(-i\theta_2)) = \frac{1}{h}i\sin\theta_2 \qquad (34)$$

are Fourier symbols of the discrete operators $-\Delta_h$, $\partial_x^h$ and $\partial_y^h$ with discrete stencils shown in (5), respectively.

*4.2 Optimal one-stage smoothing factor*

From (12) and (17), for a 2D discrete scalar operator with standard coarsening, the Fourier representation of the operator $Q_h^{2h}$ in the subspace defined in (15) is given as,

$$Q_h^{2h}\Big|_{F_{2h}(\vec{\theta})} =: \widehat{Q}_h^{2h} = \text{diag}(0,1) \in \mathbb{C}^{2\times 2}. \qquad (35)$$

Then, the smoothing factor of the discrete operator (9) is defined by [2]

$$\rho(n, D_h) = \sup_{\vec{\theta}\in\Theta_{low}} (\rho(\widehat{Q}_h^{2h}(\widehat{S}_h(\vec{\theta},\omega))^n))^{1/n}, \qquad (36)$$

i.e., the asymptotic error from the high frequency components is reduced by $n$ sweeps of relaxations, where $\widehat{S}_h(\vec{\theta},\omega)$ is the Fourier representation of the relaxation operator, $S_h(\omega)$, in the subspace shown in (15), and $\omega$ is a relaxation parameter.

We know that a good smoothing factor can be obtained by using the one-stage parameter $\omega$ in the relaxation operator $S_h(\omega)$ with the Fourier representation $\widehat{S}_h(\vec{\theta},\omega)$, that is, the optimal one-stage smoothing parameter and related smoothing factor are given by

$$\omega_{opt} = \frac{2}{2 - S_{max} - S_{min}}, \quad \rho_{opt} = \frac{S_{max} - S_{min}}{2 - S_{max} - S_{min}}, \qquad (37)$$

where $S_{max}, S_{min} \in (-1,1)$ are the maximum and minimum eigenvalues of the matrix $\widehat{Q}_h^{2h}\widehat{S}_h(\vec{\theta},1)$ with the relaxation parameter $\omega = 1$ for $\vec{\theta} \in \Theta_{low}$.

*4.3. Optimal smoothing factor for Stokes flow*

For the two-color Jacobi point relaxation $S_h^{01}$, i.e., both red and black points are swept by the Jacobi point relaxation method [1,2], applying Theorem 2, the Fourier representation of $S_h^{01}(\omega)$ with the relaxation parameter $\omega = 1$ can be given as follows by applying Theorem 2,

$$\widehat{S}_h^{01}(\vec{\theta}) = \widehat{S}_h^{01}(\vec{\theta},1) = \widehat{S}_h^1(\vec{\theta})\widehat{S}_h^0(\vec{\theta}) = \frac{1}{2}\begin{pmatrix} A_0+1 & -A_1+1 \\ -A_0+1 & A_1+1 \end{pmatrix} \cdot \frac{1}{2}\begin{pmatrix} A_0+1 & A_1-1 \\ A_0-1 & A_1+1 \end{pmatrix}, \qquad (38)$$

where

$$A_\alpha = 1 - \tilde{D}_h(\vec{\theta}^\alpha)\big/\tilde{D}_h^0(\vec{\theta}^\alpha) \qquad (39)$$

denotes the Fourier symbol of the Jacobi point relaxation for the discrete operator (9) on the grid $G_h^\beta$ shown in (14), $\tilde{D}_h^0(\vec{\theta}^\alpha)$ is the Fourier symbol of the discrete operator with the stencil $[l_{(0,0)}]_h$, and $\alpha = 0, 1$.

It is known that the smoothing factor of (6) with the distributive relaxation (28) is determined by the diagonal blocks of the transformed system (29) [2,14], which is given by

$$\rho(n, L_h) = \max\{\rho(n, -\Delta_h), \rho(n, ch^2 \Delta_h^2 - \Delta_{2h})\}. \tag{40}$$

Therefore, for the smoothing factor of (6), we only need to discuss the 2D discrete operator $-\Delta_h$ and $ch^2 \Delta_h^2 - \Delta_{2h}$. Next, the optimal smoothing factor of (6) can be investigated by using (37). For the sake of the convenient discussion, we let

$$s_1 = \sin^2 \frac{\theta_1}{2},\ s_2 = \sin^2 \frac{\theta_2}{2}. \tag{41}$$

Then, $\vec{\theta} = (\theta_1, \theta_2) \in \Theta_{low} = (-\frac{\pi}{2}, \frac{\pi}{2}]^2$ is transformed to $\vec{s} = (s_1, s_2) \in S_{low} = [0, \frac{1}{2}]^2$.

**Theorem 3** For the Poisson operator $-\Delta_h \in \mathbb{C}^{2\times 2}$, the optimal one-stage relaxation parameter and related smoothing factor of the two-color Jacobi point relaxation are $\omega_{opt} = \frac{16}{17}$, $\rho_{opt} = \frac{1}{17}$.

Proof. To investigate the two-color Jacobi point relaxation applied to the Poisson operator $-\Delta_h$ with the discrete stencil defined in (5), we substitute (16), (32) and (41) into (38) and (39), resulting in the product of (35) and (38) as follows,

$$\widehat{Q}_h^{2h} \widehat{S}_h^{01}(\vec{\theta}, 1) = \begin{pmatrix} 0 & 0 \\ \frac{1}{2}(s_1 + s_2)(1 - s_1 - s_2) & \frac{1}{2}(s_1 + s_2)(s_1 + s_2 - 1) \end{pmatrix}. \tag{42}$$

Therefore, the eigenvalues of (42) are yielded as

$$\lambda_1(s_1, s_2) = \frac{1}{2}(s_1 + s_2)(s_1 + s_2 - 1),\ \lambda_2 = 0. \tag{43}$$

So, by using (37), we obtain the optimal smoothing parameters of the two-color relaxation below,

$$S_{\max} = \max_{(s_1, s_2) \in [0, 1/2]^2} \lambda_1(s_1, s_2) \bigg|_{\vec{\theta} = (\frac{\pi}{2}, \frac{\pi}{2})} = 0,\ S_{\min} = \min_{(s_1, s_2) \in [0, 1/2]^2} \lambda_1(s_1, s_2) \bigg|_{\vec{\theta} = (0, \frac{\pi}{2})} = -\frac{1}{8}, \tag{44}$$

and

$$\omega_{opt} = \frac{2}{2 - S_{\max} - S_{\min}} = \frac{16}{17},\ \rho_{opt} = \frac{S_{\max} - S_{\min}}{2 - S_{\max} - S_{\min}} = \frac{1}{17}. \tag{45}$$

Thus the theorem is proved. □

Next, we will investigate the 2D scalar operator $ch^2 \Delta_h^2 - \Delta_{2h}$ for the two-color Jacobi point relaxation by means of the LFA. Meanwhile, the smoothing factor of the discrete system (6) can be described in the following theorem.

**Theorem 4** For the discrete operator $ch^2\Delta_h^2 - \Delta_{2h}$ with $c > 0$, the optimal one-stage relaxation parameter and related smoothing factor of the two-color Jacobi point relaxation are given as

($i$). For $c \in (0, \frac{1}{8}) \cup (\frac{1}{8}, +\infty)$,

$$\omega_{opt}(c) = \frac{3456(1-8c)^2 c^2(20c+1)^2}{\begin{bmatrix} -1 - 36c + 2736c^2 + 157248c^3 - 2115072c^4 - 2985984c^5 + 52199424c^6 + \\ (1+24c+336c^2-6912c^3+82944c^4)\sqrt{82944c^4-6912c^3+336c^2+24c+1} \end{bmatrix}},$$

$$\rho_{opt}(c) = \frac{\begin{bmatrix} -1 - 36c + 2736c^2 - 63936c^3 + 539136c^4 - 2985984c^5 - 4423680c^6 + \\ (1+24c+336c^2-6912c^3+82944c^4)\sqrt{82944c^4-6912c^3+336c^2+24c+1} \end{bmatrix}}{\begin{bmatrix} -1 - 36c + 2736c^2 + 157248c^3 - 2115072c^4 - 2985984c^5 + 52199424c^6 + \\ (1+24c+336c^2-6912c^3+82944c^4)\sqrt{82944c^4-6912c^3+336c^2+24c+1} \end{bmatrix}}.$$

($ii$). For $c = \frac{1}{8}$, $\omega_{opt}(\frac{1}{8}) = \frac{98}{217}$, $\rho_{opt}(\frac{1}{8}) = \frac{25}{217}$.

Proof. For the 2D discrete operator

$$D_h = ch^2\Delta_h^2 - \Delta_{2h}, \tag{46}$$

the Fourier symbol of (48) as follows given by applying (31)-(34),

$$\tilde{D}_h(\vec{\theta}) = \frac{1}{h^2}[4c(2 - \cos\theta_1 - \cos\theta_2)^2 + \sin^2\theta_1 + \sin^2\theta_2]. \tag{47}$$

As for the two-color Jacobi point relaxation for the discrete operator (46) with stencils (30) under consideration, we substitute (16), (32) and (41) into (38) and (39), and obtain the product of (35) and (38) below,

$$\widehat{Q}_h^{2h}\widehat{S}_h^{01}(\vec{\theta},1) = \begin{pmatrix} r_{11} & r_{12} \\ r_{21} & r_{22} \end{pmatrix}, \tag{48}$$

where $r_{11} = r_{12} = 0$, and $r_{21}$ and $r_{22}$ are expressed as

$$r_{21} = -\frac{64c(-1+s_1+s_2)\left[s_1 - s_1^2 + s_2 - s_2^2 + 4c(s_1+s_2)^2\right]}{(1+20c)^2}, \tag{49}$$

$$r_{22} = \frac{1}{(1+20c)^2} \cdot \begin{Bmatrix} 4\left[-(s_1-s_1^2+s_2-s_2^2) - 4c(-2+s_1+s_2)^2\right]\left[s_1-s_1^2+s_2-s_2^2+4c(s_1+s_2)^2\right] \\ +\left[1+20c - 2(s_1-s_1^2+s_2-s_2^2) - 8c(-2+s_1+s_2)^2\right]^2 \end{Bmatrix}. \tag{50}$$

Therefore, the eigenvalues of (48) are

$$\lambda_1 = 0, \ \lambda_2 = r_{22}. \tag{51}$$

Solving equations $\dfrac{\partial \lambda_2}{\partial s_1} = 0$ and $\dfrac{\partial \lambda_2}{\partial s_2} = 0$ on $\vec{s} = (s_1, s_2) \in [0, \dfrac{1}{2}]^2$ with $c > 0$ by using the following equation

$$\frac{\partial \lambda_2}{\partial s_1} - \frac{\partial \lambda_2}{\partial s_2} = \frac{8(s_1 - s_2)[1 + 4c(1 + 4s_1 + 4s_2)]}{(1 + 20c)^2} = 0. \tag{52}$$

which leads $s_1 = s_2$, we obtain the solution below.

$$s_1^* = s_2^* = -\frac{-1 + 36c - 480c^2 + \sqrt{1 + 24c + 336c^2 - 6912c^3 + 82944c^4}}{96c(-1 + 8c)} \tag{53}$$

with $c \in (0, \dfrac{1}{8}) \cup (\dfrac{1}{8}, +\infty)$. Then, from values of $\lambda_2(0,0)$, $\lambda_2(0, \dfrac{1}{2})$, $\lambda_2(\dfrac{1}{2}, 0)$, $\lambda_2(\dfrac{1}{2}, \dfrac{1}{2})$ and

$$\lambda_2(s_1^*, s_2^*) = \frac{\begin{bmatrix} 1 + 36c - 1008c^2 - 5184c^3 + 483840c^4 - 3649536c^5 + 20348928c^6 \\ -(1 + 24c + 336c^2 - 6912c^3 + 82944c^4)\sqrt{1 + 24c + 336c^2 - 6912c^3 + 82944c^4} \end{bmatrix}}{1728c^2(1 - 8c)^2(1 + 20c)^2},$$

for $c \in (0, \dfrac{1}{8}) \cup (\dfrac{1}{8}, +\infty)$, the maximum and minimum eigenvalues of the matrix defined in (48) are

$$S_{min} = \lambda_2(s_1^*, s_2^*) \in (-1, 0), \ S_{max} = \lambda_2(0, 0) \in (0, 1). \tag{54}$$

Furthermore, applying (37), (53) and (54), we prove Part (i) of the theorem.

When $c = \dfrac{1}{8}$, the nonzero eigenvalue of (48) is

$$\lambda_2(s_1, s_2)\bigg|_{c=\frac{1}{8}} = \frac{1}{49}\Big[1 - 40(s_1 + s_2 - s_1 s_2) + 16(s_1 - s_2)^2(s_1 + s_2) + 44(s_1^2 + s_2^2)\Big]. \tag{55}$$

Following the same way to discuss the case of $c \in (0, \dfrac{1}{8}) \cup (\dfrac{1}{8}, +\infty)$, the maximum and minimum eigenvalues of the matrix defined in (48) are

$$S_{min} = -\frac{23}{98}, \ S_{max} = \frac{1}{49}. \tag{56}$$

Due to (37), the optimal one-stage relaxation parameter and related smoothing factor are

$$\omega_{opt}(\frac{1}{8}) = \frac{2}{2 - S_{max} - S_{min}} = \frac{2}{2 - \dfrac{1}{49} + \dfrac{23}{98}} = \frac{98}{217}, \tag{57}$$

$$\rho_{opt}(\frac{1}{8}) = \frac{S_{max} - S_{min}}{2 - S_{max} - S_{min}} = \frac{\frac{1}{49} + \frac{23}{98}}{2 - \frac{1}{49} + \frac{23}{98}} = \frac{25}{217}. \tag{58}$$

Thus we prove the theorem.□

Theorem 4 results in the curve of $\rho_{opt}(c)$ that is shown in Figure 1.

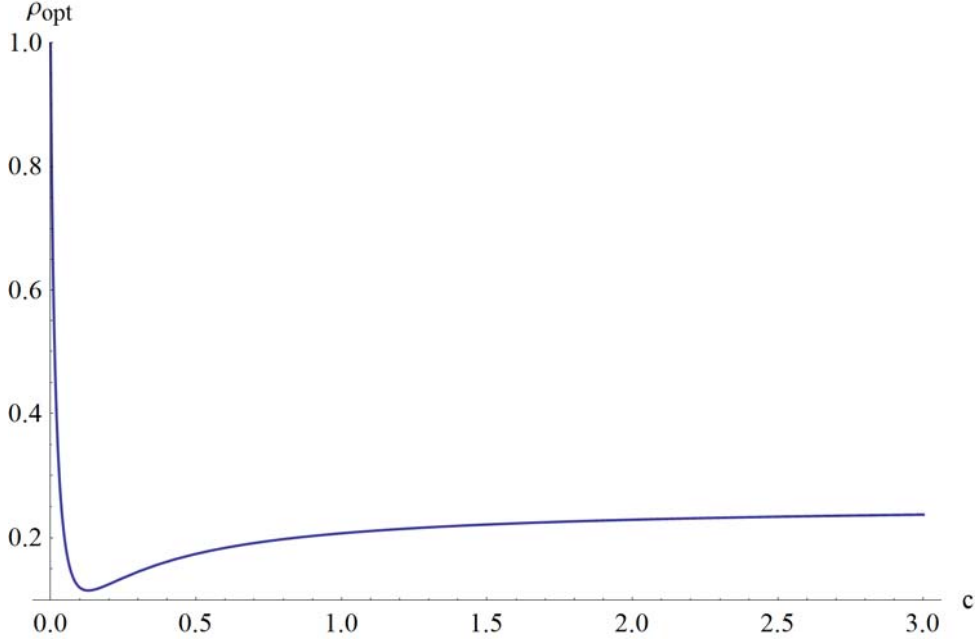

Figure 1. Curve of $\rho_{opt}(c)$

From Figure 1 and Theorem 4, we see the following relationship holds

$$\min_{c>0} \rho_{opt}(c) = \lim_{c \to \frac{1}{8}} \rho_{opt}(c) = \frac{25}{217} = \rho_{opt}(\frac{1}{8}), \tag{59}$$

and

$$\lim_{c \to +\infty} \rho_{opt}(c) = \frac{11}{43}, \quad \lim_{c \to 0} \rho_{opt}(c) = 1. \tag{60}$$

Solve $\rho_{opt}(c) = \frac{11}{43}$, yields the following parameter,

$$\frac{1}{28} < c_0 = 0.0360548 < \frac{1}{27}. \tag{61}$$

Therefore, the zone of $\rho_{opt}(c)$ is stated as

$$\frac{25}{217} \leq \rho_{opt}(c) \leq \frac{11}{43}, \ \forall \ c > \frac{1}{27}, \tag{62}$$

$$\rho_{opt}(\frac{1}{8}) < \rho_{opt}(\frac{1}{27}) \leq \rho_{opt}(c) < 1. \ \forall \ 0 < c \leq \frac{1}{27}, \tag{63}$$

And, the curve of $\omega_{opt}(c)$ is presented in Figure 2.

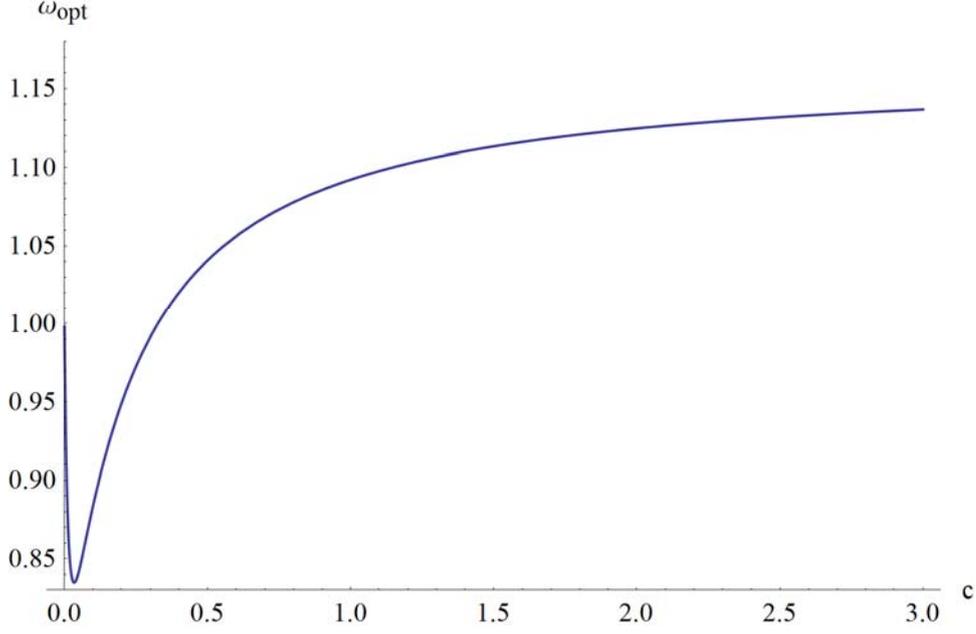

Figure 2. Curve of $\omega_{opt}(c)$

The following expression holds due to Figure 2 and Theorem 4,

$$\min_{c>0} \omega_{opt}(c) = 0.834733, \tag{64}$$

$$\lim_{c \to +\infty} \omega_{opt}(c) = \frac{50}{43}, \ \lim_{c \to 0} \omega_{opt}(c) = 1, \ \lim_{c \to \frac{1}{8}} \omega_{opt}(c) = \frac{28}{31} \neq \omega_{opt}(\frac{1}{8}). \tag{65}$$

Therefore, the zone of $\omega_{opt}(c)$ is

$$\omega_{opt}(c) \in [0.834733, \frac{50}{43}] \cup \{\frac{98}{217}\}. \tag{66}$$

According to (40), (62), (63), Theorems 3 and 4, the smoothing factor of the distributive relaxation based on the two-color Jacobi point relaxation for the discrete Stokes system is given as

$$\frac{25}{217} \leq \rho(1, L_h) = \max\{\rho(1, -\Delta_h), \rho(1, ch^2\Delta_h^2 - \Delta_{2h})\} = \rho(1, ch^2\Delta_h^2 - \Delta_{2h}) \leq \frac{11}{43}, \ \forall \ c > \frac{1}{27};$$

$$\frac{25}{217} < \rho(1, L_h) = \max\{\rho(1, -\Delta_h), \rho(1, ch^2\Delta_h^2 - \Delta_{2h})\} = \rho(1, ch^2\Delta_h^2 - \Delta_{2h}) < 1, \forall\ 0 < c \leq \frac{1}{27}.$$

## 5. Conclusion

The smoothing analysis of the distributive relaxation based on the two-color Jacobi point relaxation for solving 2D Stokes flow is investigated in this paper, where the two-color smoothing process in the frame of multigrid, the mathematical constitution on Fourier modes with various frequency components are constructed, and the Fourier representation of the two-color point relaxation is presented as well. In the meanwhile, Fourier symbols with trigonometric functions for the discrete operators and relaxations are transformed to the rational functions by (43), the smoothing process is thus simplified. Furthermore, the analytical expressions of the smoothing factor for the two-color distributive relaxation are obtained successfully, showing that the value of the smoothing factor is an upper bound of the smoothing rates, and is independent of mesh size, but depends on the parameter c in the added pressure stabilization term. All obtained results are valuable to better understand numerical performance of multigrid when a 2D Stokes flow is simulated.

## Acknowledgments

The study is supported by the National Natural Science Foundation of China (NSFC) [Grant no. 51809026] and the Doctoral Research Startup Foundation of Dali University [Grant no. KYBS201731].

## References

[1] U.Trottenberg, C.W.Oosterlee and A.schuller, Multigtid, Academic Press, New York; 2001.
[2] R.Wienands, W.Joppich, Practical Fourier Analysis for Multigrid Methods, Chapman and Hall /CRC Press; 2005.
[3] W.Briggs, V.E.Henson, S.McCormick, A Multigrid Tutorial, Society for Industrial and Applied Mathematics; 2000.
[4] W.Hackbusch, Multigrid Methods and Applications, Springer, Berlin; 1985.
[5] P.Wesseling, An Introduction to Multigrid Methods, John Wiley, Chichester, UK; 1992.
[6] K.Stuben and U.Trottenberg, Multigrid methods: Fundamental algorithms, model problem analysis andapplications, In Multigrid Methods, Volume 960 of Lectwe Notes in Mathematics, (Edited by W. Hackbusch and U. Trottenberg), pp. 1-176, Springer-Verlag, Berlin; 1982.
[7] A. Brandt, O.E.Livne, 1984 Guide to multigrid development in Multigrid Methods. Society for Industrial and Applied Mathematics; Revised edition; 2011.
[8] C.Oosterlee, F.Lorenz, Multigrid methods for the Stokes system. Comput. Sci. Eng 2006; 8(6), 34 – 43.
[9] A.Brandt, N.Dinar, Multigrid Solutions to Elliptic Llow Problems. Institute for Computer Applications in Science and Engineering, NASA Langley Research Center (1979).
[10] G.Wittum, Multigrid methods for Stokes and Navier-Stokes equations. Numerische Mathematic 1989; 54 (5), 543–563.
[11] Ming wang, Long chen, Multigrid methods for the Stokes equations using distributive Gauss–Seidel relaxations based on the least squares commutator. Journal of Scientific Computing, 2013; 56 (2), 409-431.
[12] T.Geenen, C.Vuik, G.Segal and S.MacLachlan, On iterative methods for the incompressible Stokes problem. Int. J. Numer. Methods fluids 2011; 65(10), 1180–1200.
[13] C.Bacuta, P.Vassilevski, S.Zhang, A new approach for solving Stokes systems arising from a distributive relaxation method. Numer. Methods Partial Differ. Equ 2011; 27(4), 898–914.
[14] R. Wienands, F. J. Gaspar, F. J. Lisbona and C. W. Oosterlee, An efficient multigrid solver based on distributive smoothing for poroelasticity equations. Computing 2004; 73, 99–119.
[15] Wei Liao, Boris Diskin, Yan Peng and Li-Shi Luo, Textbook-efficiency multigrid solver for three-dimensional unsteady compressible Navier-Stokes equations. Journal of Computational Physics, 227, 7160–7177 (2008).
[16] S.Takacs, All-at-once multigrid methods for optimality systems arising from optimal control problems. Ph.D. thesis, Johannes Kepler University Linz, Doctoral Program Computational Mathematics, (2012).

[17] V.Pillwein, S.Takacs, Smoothing analysis of an all-at-once multigrid approach for optimal control problems using symbolic computation. In U.Langer, P.Paule, (Eds), Numerical and Symbolic Scientific Computing: Progress and Prospects. Springer, Wien, (2011).
[18] Veronika.Pillwein, Stefan.Takacs, A local Fourier convergence analysis of a multigrid method using symbolic computation. Journal of Symbolic Computation 2014; 63, 1–20.
[19] O.E. Livne and A. Brandt, Local modes analysis of multicolor and composite relaxation Schemes, Computers and Mathematics with Applications 2004; 47, 301-317.
[20] Brown J , He Y , Maclachlan S , et al. Tuning Multigrid Methods with Robust Optimization and Local Fourier Analysis[J]. SIAM Journal on Scientific Computing, 2021, 43(1):A109-A138
[21] C.C.J. Kuo and T.F. Chan, Two-color Fourier analysis of iterative algorithms for elliptic problems with red/black ordering, SIAM J. Sci. Stat. Comput 1990; 11(4), 767-793.
[22] Irad Yavneh, On red black SOR smoothing in multigrid, SIAM J. Sci. Comput 1996;17 (l), 180-192.
[23] C.C.J. Kuo and B.C. Levy, Two-color Fourier analysis of the multigrid method with red-black Gauss-Seidel smoothing, Applied Mathematics and Computation,1989; 29, 69-87.
[24] He, Y, A generalized and unified framework of local fourier analysis using matrix-stencils. SIAM Journal on Matrix Analysis and Applications 2021; 2021(3), 42.
[25] I. Yavneh, Smoothing factors of two-color Gauss-Seidel relaxation for a class of elliptic operators, SIAM J.Num. Anal 1995; 32 (4), 1126-1138.
[26] B. Gmeinera, T. Gradl, F. Gaspar, U. Rüde, Optimization of the multigrid-convergence rate on semi-structured meshes by local Fourier analysis, Computers and Mathematics with Applications 2013; 65, 694-711.
[27] C. Rodrigo, F.J. Gaspar, F.J. Lisbona, Multicolor Fourier analysis of the multigrid method for quadratic FEM discretizations, Applied Mathematics and Computation 2012; 218,11182–11195.